\documentclass{elsart}
\usepackage{natbib}
\usepackage{amsthm, amsmath, amsfonts}

\begin{document}

\begin{frontmatter}

\title{Polynomial least squares fitting in the Bernstein basis}

\author{Ana Marco \thanksref{EM1}}
\author{, Jos\'e-Javier Mart{\'\i}nez \thanksref{EM2}}

\address{Departamento de Matem\'aticas, Universidad
de Alcal\'a,}

\address{Campus Universitario, 28871-Alcal\'a de Henares (Madrid), Spain}

\thanks[EM1]{E-mail: ana.marco@uah.es}

\thanks[EM2]{Corresponding author. E-mail: jjavier.martinez@uah.es}

\begin{abstract}

The problem of polynomial regression in which the usual monomial basis is replaced by the Bernstein basis is considered. The coefficient matrix $A$ of the overdetermined system to be solved in the least squares sense is then a rectangular Bernstein-Vandermonde matrix.

In order to use the method based on the QR decomposition of $A$, the first stage consists of computing the bidiagonal decomposition of the coefficient matrix $A$.

Starting from that bidiagonal decomposition, an algorithm for obtaining the QR decomposition of $A$ is the applied. Finally, a triangular system is solved by using the bidiagonal decomposition of the R-factor of $A$.

Some numerical experiments showing the behavior of this approach are included.

\bigskip



\end{abstract}

\end{frontmatter}

\section{Introduction}

Given $\{x_i \}_{1\leq i \leq l+1}$ pairwise distinct real points and $\{f_i \}_{1\leq i \leq l+1} \in {\bf R}$, let us consider a degree $n$ polynomial
$$
P(x)=c_0+c_1x+\ldots+c_nx^n
$$
for some $n \leq l$. Such a polynomial is a {\it least squares fit} to the data if it minimizes the sum of the squares of the deviations from the data,
$$
\sum_{i=1}^{l+1} |P(x_i)-f_i|^2.
$$
Computing the coefficients $c_j$ of that polynomial $P(x)$ is equivalent to solve, in the least squares sense, the overdetermined linear system $Ac=f$, where $A$ is the rectangular $(l+1) \times (n+1)$ Vandermonde matrix corresponding to the nodes $\{x_i \}_{1\leq i \leq l+1}$.

Taking into account that $A$ has full rank $n+1$, the problem has a unique solution given by the unique solution of the linear system
$$
A^TAc=A^Tf,
$$
the {\it normal equations}.

Since $A$ is usually an ill-conditioned matrix, it was early recognized that solving the normal equations was not an adequate method. Golub [7], following previous ideas by Householder, suggested the use of the {\it QR factorization} of $A$, which involves the solution of a linear system with the triangular matrix $R$.

Let us observe that, if $A=QR$ with $Q$ being an orthogonal matrix, then using the condition number in the spectral norm we have
$$
\kappa_2(R)=\kappa_2(A),
$$
that is, $R$ inherits the ill-conditioning of $A$ while $\kappa_2(A^TA)=\kappa_2(A)^2$.

In addition, as it was already observed by Golub in [8] (see also Section 20.1 of [9]), although the use of the orthogonal transformation avoids some of the ill effects inherent in the use of normal equations, the value $\kappa_2(A)^2$ is still relevant to some extent.

Consequently a good idea is to use, instead of the monomial basis, a polynomial basis which leads to a matrix $A$ with smaller condition number than the Vandermonde matrix.

It is frequently assumed that this happens when bases of orthogonal polynomials, such as the basis of Chebyshev polynomials, are considered. However, this fact is true when special sets of nodes are considered, but not in the case of general nodes.  A basis which leads to a matrix $A$ better conditioned than the Vandermonde matrix is the Bernstein basis of polynomials, a widely used basis in Computer Aided Geometric Designed due to the good properties that it possess (see, for instance, [2,10]). We illustrate these facts with Table $1$, where the condition numbers of Vandermonde, Chebyshev-Vandermonde and Bernstein-Vandermonde matrices are presented for the nodes considered in Example 5.1 and in Example 5.2.

\bigskip
\begin{table}[h]
\begin{center}
\begin{tabular}{|c|c|c|c|}
\hline Example & $V$ & $TV$ & $BV$  \\
\hline  5.1 & 1.7e+12 & 1.7e+12 & 2.0e+05 \\
\hline  5.2 & 2.5e+14 & 4.0e+14 & 5.3e+08 \\
\hline
\end{tabular}
\end{center}\caption{Condition numbers of the Vandermonde $V$, Chebyshev-Vandermonde $TV$ and Bernstein-Vandermonde $BV$ matrices}
\end{table}
\bigskip

Let us observe that, without lost of generality, we can consider the nodes $\{x_i \}_{1\leq i \leq l+1}$ ordered and belonging to $(0,1)$. So, we will solve the following problem:

Let $\{x_i \}_{1\leq i \leq l+1} \in (0,1)$ a set of points such that $0 <x_1 < \ldots < x_{l+1} < 1$. Our aim is to compute a polynomial
$$
P(x)=\sum_{j=0}^{n} c_{j} b_{j}^{(n)}(x)
$$
expressed in the {\it Bernstein basis} of the space $\Pi_n(x)$ of the polynomials of degree less that or equal to $n$ on the interval $[0,1]$
$$
\mathcal{B}_n=\big\{ b_j^{(n)}(x) = {n \choose j} (1 - x)^{n-j} x^j,
\qquad j = 0, \ldots, n \big\},
$$
such that $P(x)$ minimizes the sum of the squares of the deviations from the data.

This problem is equivalent to solve in the least squares sense the overdetermined linear system $Ac=f$ where now
$$
A=
\begin{pmatrix}
{n \choose 0}(1-x_1)^n & {n \choose 1}x_1(1-x_1)^{n-1}& \cdots & {n
\choose n}x_1^n\\
{n \choose 0}(1-x_2)^n & {n \choose
1}x_2(1-x_2)^{n-1}& \cdots & {n \choose n}x_2^n\\
\vdots & \vdots & \ddots & \vdots\\ {n \choose 0}(1-x_{l+1})^n & {n
\choose 1}x_{l+1}(1-x_{l+1})^{n-1}& \cdots & {n \choose n}x_{l+1}^n
\end{pmatrix} \eqno(1.1)
$$
is the $(l+1)\times(n+1)$ Bernstein-Vandermonde matrix for the Bernstein basis $\mathcal{B}_n$ and the nodes $\{x_i \}_{1\leq i \leq l+1}$,
$$f=(f_1, f_2, \ldots, f_{l+1})^T \eqno(1.2)$$ is the data vector, and
$$c=(c_1, c_2, \ldots, c_{n+1})^T \eqno(1.3)$$ is the vector containing the coefficients of the polynomial that we want to compute.

The rest of the paper is organized as follows. Neville elimination and total positivity are considered in Section 2. In Section 3, the bidiagonal factorization of a rectangular Bernstein-Vandermonde matrix is presented. The algorithm for computing the regression polynomial in Bernstein basis is given in Section 4. Finally, Section 5 is devoted to illustrate the accuracy of our algorithm by means of some numerical experiments.

\section{Basic results on Neville elimination and total positivity}

In this section we will briefly recall some basic results on Neville
elimination and total positivity which we will apply in Section 3.
Our notation follows the notation used in [4] and [5]. Given $k$,
$n \in {\bf N}$  ($1 \leq k \leq n$), $Q_{k,n}$ will denote the set
of all increasing sequences of $k$ positive integers less than or
equal to $n$.

Let $A$ be an $l \times n$ real matrix. For $k \leq l$, $m
\leq n$, and for any $\alpha \in Q_{k,l}$ and $\beta \in Q_{m,n}$,
we will denote by $A[\alpha \vert \beta ]$ the submatrix $k\times m$
of $A$ containing the rows numbered by $\alpha $ and the columns
numbered by $\beta $.

The fundamental tool for obtaining the results presented in this
paper is the {\it Neville elimination} [4,5], a procedure that
makes zeros in a matrix adding to a given row an appropriate
multiple of the previous one. We will describe the Neville elimination for a matrix $A=(a_{i,j})_{1\leq i\leq l; 1\leq j\leq n}$ where $l \geq n$. The case in which $l<n$ is analogous.

Let $A=(a_{i,j})_{1\leq i\leq l; 1\leq j\leq n}$ be a matrix where $l \geq n$. The Neville elimination of $A$ consists of $n-1$ steps
resulting in a sequence of matrices $A:=A_1\to A_2\to \ldots \to
A_n$, where $A_t= (a_{i,j}^{(t)})_{1\leq i\leq l; 1\leq j\leq n}$ has zeros below
its main diagonal in the $t-1$ first columns. The matrix $A_{t+1}$
is obtained from $A_t$ ($t=1,\ldots ,n$) by using the following
formula:
$$
a_{i,j}^{(t+1)}:= \left\{ \begin{array}{ll} a_{i,j}^{(t)}~, &
\text{if} \quad i\leq t\\
a_{i,j}^{(t)}-(a_{i,t}^{(t)}/a_{i-1,t}^{t})a_{i-1,j}^{(t)}~,~ &
\text{if} \quad i\geq t+1 ~\text{and}~ j\geq t+1\\
0~, & \text{otherwise}.
\end{array}\right.  \eqno(2.1)
$$
In this process the element
$$
p_{i,j}:=a_{i,j}^{(j)} \qquad 1\leq j\leq n, ~~ j\leq i\leq l
$$
is called {\it pivot} ($i,j$) of the Neville elimination of $A$. The
process would break down if any of the pivots $p_{i,j}$ ($1\leq j\leq n, ~~ j\leq i\leq l$) is zero. In that case we can move the corresponding rows to
the bottom and proceed with the new matrix, as described in [4].
The Neville elimination can be done without row exchanges if all the
pivots are nonzero, as it will happen in our situation. The pivots
$p_{i,i}$ are called {\it diagonal pivots}. If all the pivots
$p_{i,j}$ are nonzero, then $p_{i,1}=a_{i,1}\, \forall i$ and, by
Lemma 2.6 of [4]
$$p_{i,j}={\det A[i-j+1,\ldots ,i\vert 1,\ldots ,j]\over \det
A[i-j+1,\ldots ,i-1\vert 1,\ldots ,j-1]}\qquad 1<j\leq n, ~ j\leq i\leq l.
\eqno(2.2)$$
The element
$$
m_{i,j}=\frac{p_{i,j}}{p_{i-1,j}} \qquad 1\leq j\leq n, ~~ j< i\leq
l  \eqno(2.3)
$$
is called {\it multiplier} of the Neville elimination of $A$. The
matrix $U:=A_n$ is upper triangular and has the diagonal pivots in
its main diagonal.

The {\it complete Neville elimination} of a matrix $A$ consists on
performing the Neville elimination of $A$ for obtaining $U$ and then
continue with the Neville elimination of $U^T$. The pivot
(respectively, multiplier) $(i,j)$ of the complete Neville
elimination of $A$ is the pivot (respectively, multiplier) $(j,i)$
of the Neville elimination of $U^T$, if $j\ge i$. When no row
exchanges are needed in the Neville elimination of $A$ and $U^T$, we
say that the complete Neville elimination of $A$ can be done without
row and column exchanges, and in this case the multipliers of the
complete Neville elimination of $A$ are the multipliers of the
Neville elimination of $A$ if $i\ge j$ and the multipliers of the
Neville elimination of $A^T$ if $j\ge i$.

A matrix is called {\it totally positive} (respectively, {\it
strictly totally positive}) if all its minors are nonnegative
(respectively, positive). The Neville elimination characterizes the
strictly totally positive matrices as follows [4]:

{\bf Theorem 2.1.} A matrix is strictly totally positive if and only
if its complete Neville elimination can be performed without row and
column exchanges, the multipliers of the Neville elimination of $A$
and $A^T$ are positive, and the diagonal pivots of the Neville
elimination of $A$ are positive.

It is well known [2] that the Bernstein-Vandermonde matrix is a
strictly totally positive matrix when the nodes
satisfy $0<x_1<x_2<\ldots<x_{l+1}<1$,  but this result will also be shown to
be a consequence of our Theorem 3.2.

\section{Bidiagonal factorization of $A$}

In this section we consider the bidiagonal factorization of the Bernstein-Vandermonde matrix $A$ of (1.1).

Let us observe that when $l=n$ this matrix $A$ is the coefficient matrix of the linear system
associated with a Lagrange interpolation problem in the Bernstein basis $\mathcal{B}_n$ whose interpolation nodes are $\{x_i: ~ i=1, \ldots, n+1\}$. A fast and accurate algorithm for solving this linear system, and therefore the corresponding Lagrange interpolation problem in the Bernstein basis can be found in [13]. A good introduction to the interpolation theory can be seen in [3].

The following two results will be the key to construct our algorithm.

\medskip
{\bf Proposition 3.1.} (See [13]) Let $A$ be the square Bernstein-Vandermonde matrix of order $n+1$ for the Bernstein basis $\mathcal{B}_n$  and the nodes $x_1, x_2, \ldots, x_{n+1}$.
We have:
$$
\det A = {n \choose 0}{n \choose 1} \cdots {n \choose n}
\prod_{1\leq i < j \leq n+1}(x_j-x_i).
$$

\medskip
{\bf Theorem 3.2.} Let $A=(a_{i,j})_{1\le i \leq l+1; 1\leq j\le n+1}$ be a
Bernstein-Vandermonde matrix for the Bernstein basis $\mathcal{B}_n$ whose nodes satisfy $0 < x_1 < x_2
<\ldots < x_l <x_{l+1} <1$. Then $A$ admits a factorization in
the form
$$A=F_{l}F_{l-1} \cdots F_1 D G_1 \cdots G_{n-1}G_{n}  \eqno(3.1)$$
where $G_j$ are $(n+1)\times (n+1)$ upper triangular bidiagonal matrices ($j=1,\ldots,n$), $F_i$ are $(l+1)\times(l+1)$
lower triangular bidiagonal matrices ($i=1,\ldots,l$), and $D$ is a $(l+1) \times (n+1)$
diagonal matrix.

{\bf Proof.} The matrix $A$ is strictly totally positive (see [2])
and therefore, by Theorem 2.1, the complete Neville elimination of
$A$ can be performed without row and column exchanges providing the
following factorization of $A$ (see [6]):
$$A=F_{l}F_{l-1} \cdots F_1 D G_1 \cdots G_{n-1}G_{n},  $$
where $F_i$ ($1\le i\le l$) are $(l+1)\times(l+1)$  bidiagonal matrices of the form
$$F_i=\begin{pmatrix}
1 & & & & & & & \\
0 & 1 & & & & & & \\
& \ddots & \ddots & & & & & \\
& & 0 & 1 & & & & \\
& & & m_{i+1,1} & 1 & & & \\
& & & & m_{i+2,2} & 1 & & \\
& & & & & \ddots & \ddots & \\
& & & & & & m_{l,l-i} & 1
\end{pmatrix},  \eqno(3.2)
$$
$G^T_i$ ($1\le i\le n$) are $(n+1)\times (n+1)$ bidiagonal matrices of the form
$$G_i^T=\begin{pmatrix}
1 & & & & & & & \\
0 & 1 & & & & & & \\
& \ddots & \ddots & & & & & \\
& & 0 & 1 & & & & \\
& & & \widetilde m_{i+1,1} & 1 & & & \\
& & & & \widetilde m_{i+2,2} & 1 & & \\
& & & & & \ddots & \ddots & \\
& & & & & & \widetilde m_{n,n-i} & 1
\end{pmatrix},  \eqno(3.3)
$$
and $D$ is the $(l+1)\times (n+1)$ diagonal matrix whose $i$th ($1\le i\le n+1$)
diagonal entry is the diagonal pivot $p_{i,i}=a_{i,i}^{(i)}$ of the
Neville elimination of $A$:
$$
D =(d_{i,j})_{1\le i \leq l+1; 1\leq j\le n+1}=\text{diag}\{p_{1,1},p_{2,2},\ldots,p_{n+1,n+1}\}. \eqno(3.4)
$$

Taking into account that the minors of $A$ with $j$ initial
consecutive columns and $j$ consecutive rows starting with row $i$
are
$$
\begin{aligned}
& \det A [i,\ldots, i+j-1 \vert 1, \ldots, j] ={n \choose 0} {n
\choose 1} \cdots {n \choose j-1}\\ & (1-x_i)^{n-j+1}
(1-x_{i+1})^{n-j+1} \cdots (1-x_{i+j-1})^{n-j+1} \prod_{i \leq k <
h \leq i+j-1}(x_h-x_k),
\end{aligned}
$$
a result that follows from the properties of the determinants and
Proposition 3.1, and that $m_{i,j}$ are the multipliers of the Neville
elimination of $A$, we obtain that
$$
m_{i,j}=\frac{p_{i,j}}{p_{i-1,j}}=\frac{ (1-x_i)^{n-j+1} (1-x_{i-j})
\prod_{k=1}^{j-1}(x_i-x_{i-k}) }{ (1-x_{i-1})^{n-j+2}
\prod_{k=2}^j(x_{i-1}-x_{i-k}) },  \eqno(3.5)
$$
where $j=1,\ldots,n+1$ and $i=j+1, \dots, l+1$.

As for the minors of $A^T$ with $j$ initial consecutive columns and
$j$ consecutive rows starting with row $i$, they are:
$$
\begin{aligned}
& \det A^T[i,\ldots, i+j-1 \vert 1, \ldots, j] ={n \choose i-1} {n
\choose i} \cdots {n \choose i+j-2}x_1^{i-1}x_2^{i-1} \cdots
x_j^{i-1}\\
& (1-x_1)^{n-i-j+2}(1-x_2)^{n-i-j+2} \cdots
(1-x_j)^{n-i-j+2}\prod_{1 \leq k < h \leq j}(x_h-x_k).
\end{aligned}
$$
This expression also follows from the properties of the determinants
and Proposition 3.1. Since the entries $\widetilde m_{i,j}$ are the
multipliers of the Neville elimination of $A^T$, using the previous
expression for the minors of $A^T$ with initial consecutive columns
and consecutive rows, it is obtained that
$$
\widetilde m_{i,j}=\frac{(n-i+2)\cdot x_j}{(i-1)(1-x_j)} \qquad
\quad j=1,\dots,n;
 ~ i=j+1,\ldots,n+1. \eqno(3.6)
$$
Finally, the $i$th diagonal element of $D$
$$
p_{i,i}=\frac{ {n \choose i-1}(1-x_i)^{n-i+1} \prod_{k<i}(x_i-x_k)
}{ \prod_{k=1}^{i-1}(1-x_k) } ~, \qquad i=1,\ldots,n+1 \eqno(3.7)
$$
is obtained by using the expression for the minors of $A$ with
initial consecutive columns and initial consecutive rows. $\Box$

Moreover, by using the same arguments of [14] it can be seen that
this factorization is unique among factorizations of this type, that
is to say, factorizations in which the matrices involved have the
properties shown by formulae (3.2), (3.3) and (3.4).

\medskip
{\bf Remark 3.3.} The formulae obtained in the proof of Theorem
3.3 for the minors of $A$ with $j$ initial consecutive columns and
$j$ consecutive rows, and for the minors of $A^T$ with $j$ initial
consecutive columns and $j$ consecutive rows show that they are not
zero and so, the complete Neville elimination of $A$ can be
performed without row and column exchanges. Looking at equations
(3.5)-(3.7) is easily seen that $m_{i,j}$, $\widetilde m_{i,j}$ and
$p_{i,i}$ are positive. Therefore, taking into account Theorem 2.1,
this confirms that the matrix $A$ is strictly totally positive.

\medskip
{\bf Remark 3.4.} In the square case, the matrices $F_i$ ($i=1,\ldots, l$) and the matrices $G_j$ ($j=1,\ldots,n$) are not the same bidiagonal matrices that appear in the bidiagonal factorization of $A^{-1}$ presented in [13], nor their inverses. The multipliers of the Neville elimination of $A$ and $A^T$ give us the bidiagonal factorization of $A$ and $A^{-1}$, but obtaining the bidiagonal factorization of $A$ from the bidiagonal factorization of $A^{-1}$ (or vice versa) is not straightforward [6]. The structure of the bidiagonal matrices that appear in both factorizations is not preserved by the inversion, that is, in general, $F_i^{-1}$ ($i=1,\ldots,l$) and $G_j^{-1}$ ($j=1,\ldots,n$) are not bidiagonal matrices. See [6] for a more detailed explanation.

\section{The algorithm}

In this section we present an accurate and efficient algorithm for solving the problem of polynomial regression in Bernstein basis we have presented in Section 1. As we introduced there, our algorithm is based on the solution of the least squares problem $min_c\parallel Ac-f \parallel$, where $A$, $f$ and $c$ are given by (1.1), (1.2) and (1.3), respectively. Taking into account that $A$ is a strictly totally positive matrix, it is full rank, and the method based on the QR decomposition is the most adequate for solving this least squares problem [1].

The following result (see Section 1.3.1 in [1]) will be essential in the construction of our algorithm.

\medskip
{\bf Theorem 4.1.} Let $Ac=f$ a linear system where $A \in {\bf R}^{(l+1) \times (n+1)}$, $l \geq n$, $c \in {\bf R}^{n+1}$ and $f \in {\bf R}^{l+1}$. Assume that $rank(A)=n+1$, and let the QR decomposition of $A$ given by
$$
A=Q \begin{pmatrix} R\\ 0 \end{pmatrix},
$$
where $Q \in {\bf R}^{(l+1) \times (l+1)}$ is an orthogonal matrix and $R \in {\bf R}^{(n+1) \times (n+1)}$ is an upper triangular matrix with nonnegative diagonal entries.

The solution of the least squares problem $min_c \parallel Ac-f \parallel_2$ is obtained from
$$
\begin{pmatrix} d_1 \\ d_2\end{pmatrix} = Q^T f, \quad Rc=d_1, \quad r=Q \begin{pmatrix}0 \\ d_2 \end{pmatrix},
$$
where $d_1 \in {\bf R}^{n+1}$, $d_2\in {\bf R}^{l-n}$ and $r=f-Ac$. In particular $\parallel r \parallel_2 = \parallel d_2 \parallel_2$.

\medskip

An accurate and efficient algorithm for computing the QR decomposition of a strictly totally positive matrix $A$ is presented in [12]. This algorithm is called {\tt TNQR} and can be obtained from [11]. Given the bidiagonal factorization of $A$, {\tt TNQR} computes the matrix $Q$ and the bidiagonal factorization of the matrix $R$. Let us point out here that if $A$ is strictly totally positive, then $R$ is strictly totally positive. {\tt TNQR} is based on Givens rotations, has a computational cost of $O(l^2n)$ arithmetic operations if the matrix $Q$ is required, and its high relative accuracy comes from the avoidance of subtractive cancellation.

\medskip

A fast and accurate algorithm for computing the bidiagonal factorization of the rectangular Bernstein-Vandermonde matrix that appears in our problem of polynomial regression in the Bernstein basis can be developed by using the expressions (3.5), (3.6) and (3.7) for the computation of the multipliers $m_{i,j}$ and  $\widetilde m_{i,j}$, and the diagonal pivots $p_{i,i}$ of its Neville elimination. The algorithm is an extension to the rectangular case of the one presented in [13] for the square Bernstein-Vandermonde matrices. Given the nodes $\{x_i \}_{1\leq i \leq l+1} \in (0,1)$ and the degree $n$ of the Bernstein basis, it returns a matrix $M \in {\bf R}^{(l+1)\times(n+1)}$ such that
$$
\begin{aligned}
M_{i,i} & =p_{i,i} \quad ~~i=1,\ldots,n+1,\\
M_{i,j} & =m_{i,j} \quad  j=1,\ldots, n+1; ~i=j+1,\ldots, l+1,\\
M_{i,j} & =\widetilde m_{j,i} \quad i=1,\ldots,n; ~j=i+1,\ldots, n+1.
\end{aligned}
$$
The algorithm, that we call it {\tt TNBDBV}, has a computational cost of $O(ln)$ arithmetic operations, and high relative accuracy because it only involves arithmetic operations that avoid subtractive cancellation (see [13] for the details). The implementation in \textsc{Matlab} of the algorithm in the square case can be taken from [11].

\medskip

In this way, the algorithm for solving the least squares problem $min_c \parallel Ac-f \parallel$ corresponding to our polynomial regression problem will be:

INPUT: The nodes $\{x_i \}_{1\leq i \leq l+1} \in (0,1)$, the data vector $f$ and the degree $n$ of the Bernstein basis.

OUTPUT: A vector $c=(c_j)_{1\leq i \leq n+1}$ containing the coefficients of the polynomial $P(x)$ in the Bernstein basis $\mathcal{B}_n$ and the minimum residual $r$.

\begin{itemize}

\item[-] {\it Step 1}: Computation of the bidiagonal factorization of $A$ by means of {\tt TNBDBV}.

\item[-] {\it Step 2}: Given the matrix $M$ obtained in Step 1, computation of the QR decomposition of $A$ by using {\tt TNQR}.

\item[-] {\it Step 3}: Computation of
$$d=\begin{pmatrix}d_1 \\d_2 \end{pmatrix}=Q^Tf.$$

\item[-] {\it Step 4}: Solution of the upper triangular system $Rc=d_1$.

\item[-] {\it Step 5}: Computation of
$$
r=Q \begin{pmatrix}0 \\ d_2 \end{pmatrix}.
$$

\end{itemize}

Step 3 and Step 5 are carried out by using the standard matrix multiplication command of \textsc{Matlab}. As for Step 4, it is done by means of the algorithm {\tt TNSolve} of P. Koev [11]. Given the bidiagonal factorization of a totally positive matrix $A$, {\tt TNSolve} solves a linear system whose coefficient matrix is $A$ by using backward substitution.

Let us observe that $A$ is not constructed, although we are also computing the residual $r=f-Ac$.

\section{Numerical experiments and final remarks}

 Two numerical experiments illustrating the good properties of our algorithm are reported in this section. We solve the least squares problem $min_c\parallel Ac-f \parallel$ corresponding to the computation of the regression polynomial in exact arithmetic by means of the command {\tt leastsqrs} of {\it Maple 10} and we denote this solution by  $c_e$. We also compute the minimum residual $r_e$ in exact arithmetic by using {\it Maple 10}. We use $c_e$ and $r_e$ for comparing the accuracy of the results obtained in \textsc{Matlab} by means of:

\begin{enumerate}
\item The algorithm presented in Section 4.

\item The command $A\backslash f$ of \textsc{Matlab}.
\end{enumerate}

The relative errors obtained when using the approaches (1) and (2) for computing the coefficients of the regression polynomial in the experiments described in this section ($ec_1$ and $ec_2$, respectively) are included in the first and in the third column of Table 2. The relative errors corresponding to the computation of the minimum residual by using the approaches (1) and (2) ($er_1$ and $er_2$, respectively) are presented in the second and in the fourth column of Table 2.

We compute the relative error of a solution $c$ of the least squares problem $min_c\parallel Ac-f \parallel$ by means of the formula
$$
ec=\frac{\parallel c-c_e \parallel_2}{\parallel c_e \parallel_2}.
$$
The relative error of a minimum residual $r$ is computed by means of
$$
er=\frac{\parallel r-r_e \parallel_2}{\parallel r_e \parallel_2}.
$$

{\bf Example 5.1.} Let $\mathcal{B}_{15}$ the Bernstein basis of the space of polynomials with degree less than or equal to $15$ in $[0,1]$. We will compute the polynomial
$$
P(x)=\sum_{j=0}^{15} c_{j} b_{j}^{(n)}(x)
$$
that minimizes
$$
\sum_{i=1}^{21} |P(x_i)-f_i|^2,
$$
where
$$
\{x_i\}_{1 \leq i \leq 21} = \Big\{ \frac{i}{22} \Big\}_{1 \leq i \leq 21},
$$
and
$$
f=(3,4,0,-2,5,0,1,9,-3,7,-1,0,2,2,-4,-2,3,8,-6,4,1)^T.
$$
Let us observe that, the condition number of the Bernstein-Vandermonde matrix $A$ of the least squares problem corresponding to the regression polynomial we are interested in computing is $\kappa_2(A)=2.0e+05$.

The following example shows how the algorithm we have presented in this paper keeps the accuracy when the condition number of the Bernstein-Vander-\break monde matrix involved in the regression problem increases, while the accuracy of the general approach (2) which does not exploit the structure of this matrix goes down.

\medskip
{\bf Example 5.2.} We consider a regression problem such that the Bernstein basis $\mathcal{B}_{15}$ and the data vector $f$ are the same as in Example 5.1. The points $\{x_i\}_{1 \leq i \leq 21}$ are now:
$$
\Big\{ \frac{1}{22}, \frac{1}{20}, \frac{1}{18}, \frac{1}{16}, \frac{1}{14}, \frac{1}{12}, \frac{1}{10}, \frac{1}{8}, \frac{1}{6}, \frac{1}{4}, \frac{1}{2}, \frac{23}{42}, \frac{21}{38}, \frac{19}{34}, \frac{17}{30}, \frac{15}{26}, \frac{13}{22}, \frac{11}{18}, \frac{9}{14}, \frac{7}{10}, \frac{5}{6} \Big\}.$$

The condition number of the Bernstein-Vandermonde matrix $A$ involved in this experiment is $\kappa_2(A)=5.3e+08$.

\begin{table}[h]
\begin{center}
\begin{tabular}{|c|c|c|c|c|}
\hline Example & $ec_1$ & $er_1$ & $ec_2$ & $er_2$ \\
\hline  5.1 & 1.4e-15 & 1.3e-15 & 9.0e-12 & 5.2e-12\\
\hline  5.2 & 2.0e-15 & 2.3e-15 & 1.9e-09 & 1.4e-08\\
\hline
\end{tabular}
\end{center}\caption{Relative errors in Example 5.1 and Example 5.2}
\end{table}

\medskip

{\bf Remark 5.3.} The accuracy of our algorithm is obtained by exploiting the structure of the Bernstein-Vandermonde matrix. Every step of our algorithm, except the ones in which the standard matrix multiplication command of \textsc{Matlab} is used, are developed with high relative accuracy because only arithmetic operations that avoid subtractive cancellation are involved [12,13].

\medskip
{\bf Remark 5.4.} Our algorithm  has the same computational cost ($O(l^2n)$ arithmetic operations) as the conventional algorithms that solve the least squares problem by means of the QR decomposition ignoring the structure of the matrix, when $Q$ is explicitly required (see Section 2.4.1 of [1]).

%Although, of course, this problem can be solved by using standard algorithms for least squares problems, the fact that Bernstein-Vandermonde matrices %are ill conditioned \cite{MM07} makes these algorithms not to give accurate solutions.

\end{document}